\documentclass[11pt,a4paper]{amsart} 

\usepackage{latexsym}

\usepackage[a4paper , lmargin = {2.5cm} , rmargin = {2.5cm} , tmargin = {2.5cm} , bmargin = {2.5cm} ]{geometry} 

\usepackage{geometry}
\usepackage{pinlabel}         
\usepackage{graphicx}
\usepackage{amssymb, amsthm}
\usepackage{epstopdf}
\usepackage{romannum}
\usepackage{tikz}
\usepackage{subfig}
\usepackage{MnSymbol}
\usepackage{ulem}
\usepackage{cancel}

  	\newcommand{\Z}{\ensuremath{\mathbb{Z}}}

   	\def\Aut{{\rm{Aut}}$(F_n)$}
        \def\Autn{{\rm{Aut}}(F_n)}
        
   	\def\SAut{{\rm{SAut}}$(F_n)$}
        \def\SAutn{{\rm{SAut}}(F_n)}
	\def\GL{{\rm{GL}}}
	\def\SL{{\rm{SL}}}
        
	\def\isom{{\rm{Isom}}(X)}

	\def\Fix{{\rm{Fix}}}

	\def\CAT{{\rm{CAT}}$(0)$}

	\def\F{{\rm{F}}}

\theoremstyle{plain}
\newtheorem*{NewTheorem}{Theorem}	
\newtheorem*{NewTheoremA}{Theorem A}
\newtheorem*{NewTheoremB}{Theorem B}

\newtheorem*{CellularProposition}{Fixed Point Criterion}

\newtheorem*{Question}{Question}

\newtheorem{theorem}{Theorem}[section]

\newtheorem{example}[theorem]{Example}
\newtheorem{proposition}[theorem]{Proposition}
\newtheorem{Definition}[theorem]{Definition}

\newtheorem*{HellyCubeTheorem}{Helly's Theorem for  $\mathbf{CAT(0)}$ cube complexes}

\newenvironment{pf}{\par\medskip\noindent\textit{Proof.}~}{\hfill $\square$\par\medskip}

\title{A condition that prevents groups from acting fixed point free  on  cube complexes}
\author{Olga Varghese}
\thanks{Research partially supported by SFB 878.}
\date{\today}
\address{Olga Varghese\\
Department of Mathematics\\
M\"unster University\\ 
Einsteinstra\ss e 62\\
48149 M\"unster (Germany)}
\email{olga.varghese@uni-muenster.de}

\begin{document}

\pagenumbering{arabic}
\begin{abstract}

We describe a group theoretic condition which ensures that any strongly simplicial action of a group satisfying this condition on a \CAT\ cube complex has a global fixed point. In particular, we show that this fixed point criterion is satisfied by \Aut, the automorphism group of a free group of rank $n$. For \SAut, the unique subgroup of index two in \Aut, we obtain a similar result.



\end{abstract}

\maketitle

\section{Introduction}

In the mathematical world, this article is located in the area of geometric group theory, a field at the intersection of algebra, geometry and topology. Geometric group theory studies the interaction between algebraic and geometric properties of groups. One~is interested in understanding on which 'nice' geometric spaces a given  group can act in a reasonable way and how geometric properties of these spaces are reflected in the algebraic structure of the group.
Here, the spaces will  be \CAT\ cube complexes,  while the groups will be \Aut\ and \SAut. 
The questions we shall investigate are concerned with fixed point properties of these groups.

More precisely, let $\mathbb{Z}^n$ be the free abelian group and $F_{n}$ the free group of rank $n$. One goal for a group theorist is to unterstand the structure of their automorphism groups, 
${\rm GL}_n(\mathbb{Z})$, resp. \Aut.
The abelianization map $F_{n}\twoheadrightarrow \Z^{n}$ gives a natural epimorphism 
$\Autn\twoheadrightarrow\GL_{n}(\Z).$
The special automorphism group of $F_n$, which we will denote by \SAut, is defined as the preimage of ${\rm SL}_{n}(\Z)$ under this map. 
Much of the work on \Aut\ and \SAut\ is motivated by the idea that ${\rm GL}_n(\mathbb{Z})$ and \Aut, resp. ${\rm SL}_n(\Z)$ and \SAut,  should have many properties in common. 
Here we follow this idea and present analogies between these groups with respect to fixed point properties.

The starting point for our investigation is the study of group actions on simplicial trees which was initiated by Serre, \cite{Serre1}, \cite{Serre2}. Let $\mathcal{A}$ be the class of simplicial trees. A group $G$ is said to have property $\F\mathcal{A}$ if any simplicial action of $G$ on any member of $\mathcal{A}$ has a fixed point. It was proven by Serre that ${\rm GL}_n(\Z)$ and~${\rm SL}_n(\Z)$ have property $\F\mathcal{A}$ for $n\geq 3$. Regarding \Aut\ and \SAut, Bogopolski was the first to prove that these groups also have property $\F\mathcal{A}$, see \cite{BogopolskiFA}.

In this paper, we investigate some higher dimensional analog of this property. More precisely, we say that a group $G$ has property $\F\mathcal{C}$ (resp. $\F\mathcal{C}_{*}$) if it satisfies the following property:  Every simplicial action of $G$ on a \CAT\ cube complex (resp. on a finite dimensional \CAT\ cube complex) has a fixed point.

It was noticed by Farb in \cite{Farb} that ${\rm GL}_n(\Z)$ and~${\rm SL}_n(\Z)$ have property $\F\mathcal{C}_{*}$ for $n\geq 3$. Concerning \Aut\ Bridson and Vogtmann formulated the following question.
\begin{Question}(\cite[14]{BridsonFAd})
If $n\geq 4$, can \Aut\ act without a global fixed point on a finite dimensional \CAT\ cube complex?
\end{Question} 
Under the additional assumption that the action is strongly simplicial, meaning that the stabilizer group of any cube fixes that cube pointwise, we show the answer to this question to be negative.

We say that a group $G$ has property $\F^{s}\mathcal{C}$ if any strongly simplicial action  of $G$ on any \CAT\ cube complex  has a fixed point. Note that property $\F\mathcal{C}$ trivially implies $\F^{s}\mathcal{C}$. The converse is not true, see Example \ref{example}. 

We give a group theoretic condition that implies that a group satisfying this condition has property $\F^{s}\mathcal{C}$.
\begin{CellularProposition}
Let $G$ be a group and $Y$ a finite generating set of $G$. If each pair of elements in $Y$ generates a finite subgroup, then $G$ has property $\F^{s}\mathcal{C}$.
\end{CellularProposition}
Our proof of this Fixed Point Criterion is based on a suitable version of Helly's Theorem, one important result of convexity theory. 

There exist several variations of Helly's Theorem in the literature, e.g. for finite families of convex open, resp. closed, subsets of a \CAT\ metric space, see \cite[2]{Debrunner},  \cite[3.2]{Farb} and \cite[5.3]{Kleiner}.  Indeed, it was Farb who discovered the connection between Helly's Theorem and the combinatorics of generating sets. 

We show that the Fixed Point Criterion is satisfied by the groups \Aut\ and~\SAut. We hence obtain 
\begin{NewTheoremA}
$ $
\begin{enumerate}
\item[$(i)$] For $n\geq 3$ the group \Aut\ has property $\F^{s}\mathcal{C}$. 
\item[$(ii)$] For $n\geq 4$ the group \SAut\ has property $\F^{s}\mathcal{C}$. 
\end{enumerate}
\end{NewTheoremA}
We use these results to show that actions of \SAut\ on certain low dimensional \CAT\ cube complexes are automatically trivial. 

\begin{NewTheoremB}
Let $n\geq 4$ and $X$ be a \CAT\ cube complex such that every vertex has at most $m$ neighbours. Let $\Phi:\SAutn\rightarrow{\rm Isom}(X)$
be a strongly simplicial action. If $m<n$, then $\Phi$ is trivial.
\end{NewTheoremB}

\subsection*{Remark on Kazhdan's property (T) }

Another fixed point property of interest is Kazhdan's property (T). Groups with property (T) play an important role in many areas of mathematics and even in computer science. 

It is a fundamental open question in geometric group theory whether the groups \Aut\ and \SAut\ for $n\geq 6$ have property (T) or not. Recently, it was proven by Kaluba, Nowak and Ozawa in \cite{Kaluba} that 
${\rm Aut}(F_5)$ and ${\rm SAut}(F_5)$ have property (T). Property (T) is related to the ones studied here by the following result of Niblo and Reeves.
\begin{NewTheorem}(\cite{Niblo})
If $G$ is a finitely generated group satisfying Kazhdan's property~(T), then $G$ has property $\F\mathcal{C}_*$. 
\end{NewTheorem}

A first step in answering this open question would be to show that the groups \Aut\ and \SAut\ have property $\F\mathcal{C}_*$.  

\section{ \CAT\ cube complexes}

The purpose of this subsection is to introduce \CAT\ cube complexes. A detailed description of \CAT\ spaces and their geometry can be found in  \cite{Haefliger}. Roughly speaking, a cube complex is a space which one obtains by taking a union of unit cubes of possibly different dimensions and gluing them along isometric faces. 

Let us give the formal definition of a cube complex. We denote by ${\rm I}^{d}$ the unit $d$-dimensional cube $[0,1]^{d}$. By convention, ${\rm I}^{0}$ is a point. A face of ${\rm I}^{d}$ is a subset $F$ of ${\rm I}^{d}$, that is a product $F^1\times F^2\times\ldots \times F^d$ where each $F^{i}$ is  either $\left\{0\right\}, \left\{1\right\}$ or $[0,1]$. Let $C_1, C_2$ be two cubes with faces $F_1\subset C_1, F_2\subset C_2$. A glueing of $C_1$ and $C_2$ along $F_1, F_2$ is a bijective isometry $\psi_{C_1,C_2}: F_1\rightarrow F_2$.

\begin{Definition}
Let $\mathcal{C}$ be a family of cubes and $\mathcal{F}$ be a family of glueings of elements of $\mathcal{C}$ with the properties that no cube is glued with itself and that for all
cubes $C_1$ and $C_2$ there is at most one gluing $\psi_{C_1, C_2}$.

A cube complex $K$ is the quotient of the disjoint union of cubes $X= \bigsqcup\mathcal{C}$ by the glueing equivalence relation that is generated by
\[
\left\{x\sim\psi_{C_1, C_2}(x)\mid \psi_{C_1, C_2}\in\mathcal{F}, x\in{\rm domain}(\psi_{C_1, C_2})\right\}.
\]

A subset $L$ of $K$ is a cube subcomplex of $K$ if there exist a subset $\mathcal{C}'$  of $\mathcal{C}$ such that
$L$ is equal to the quotient of a disjoint union of cubes $Y= \bigsqcup\mathcal{C}'$ by a glueing equivalence relation which is generated by
\[
\{x\sim\psi_{C_1, C_2}(x)\mid C_1, C_2\in\mathcal{C}', \psi_{C_1, C_2}\in\mathcal{F}, x\in{\rm domain}(\psi_{C_1, C_2})\}.
\]
\end{Definition}
The cube complex is \CAT\ if the cube complex with the length metric is a \CAT\ space. We note that the \CAT\ inequality condition for a cube complex can be expressed by a combinatorial condition on the cells, see  \cite[\Romannum{2}.5]{Haefliger}.

For example, the  $d$-dimensional Euclidean space $\mathbb{R}^d$ is a \CAT\ cube complex in the obvious way with $\mathbb{Z}^d$ as a set of vertices.

The following version of the Bruhat-Tits Fixed Point Theorem valid for non necessarily complete \CAT\ cube complexes was proven by  \textsc{Gerasimov} in \cite{Gerasimov}, (see \cite[5.18]{Cornulier} for details). 
\begin{proposition}
\label{boundedOrbit}
 Let $G$ be a finite group acting simplicially on a \CAT\ cube complex.  Then $G$ has a global fixed point.
\end{proposition}

\section{Helly's Theorem for cube complexes}
In this subsection we prove a version of an important result of convexity theory, Helly's Theorem for the class of \CAT\ cube complexes. 

We start by giving the following definition and result about median graphs.
\begin{Definition}
Let $\Gamma$ be a graph. The interval ${\rm I}(u,v)$ between two vertices $u$ and~$v$ consists of all vertices on a shortest paths between $u$ and $v$, i.e.
\[
{\rm I}(u,v):=\left\{x\in V\mid d(u,x)+d(x,v)=d(u,v)\right\}.
\]
A graph $\Gamma$ is called {\bf median} if for each triple $x, y, z$ of vertices the interval intersection consists of exactly one vertex, denoted by $m(x, y, z)$, i.e. 
\[
{\rm I}(x,y)\cap {\rm I}(y, z)\cap {\rm I}(x, z)=\left\{m(x,y,z)\right\}.
\]
\end{Definition}

The relation between \CAT\ cube complexes and median graphs is as follows.
\begin{proposition}(\cite[\S 10]{Roller})
\label{MG}
Let $X$ be a \CAT\ cube complex and $X^{(1)}$ be the $1$-skeleton of~$X$. Then $X^{(1)}$ is a median graph.
\end{proposition}

With the help of above proposition we can now prove a suitable version of Helly's Theorem. 

\begin{HellyCubeTheorem}(\cite[2.2]{Roller})
\label{HellyCube}
Let $X$ be a \CAT\ cube complex and $\mathcal{S}$ a finite family of non-empty  \CAT\ subcomplexes of $X$. If the intersection of each two elements of $\mathcal{S}$ is non-empty, then $\bigcap\mathcal{S}$ is non-empty.
\end{HellyCubeTheorem}

\begin{pf}
We argue by induction on $m:=|\mathcal{S}|$. For $m=2$ there is nothing to prove. Let $m\geq 3$ and 
\[
Y:=X_3\cap X_4\cap\ldots\cap X_{m}.
\]
By induction hypothesis $X_1\cap Y$, $X_2\cap Y$ and $X_1\cap X_2$ are non-empty.
Choose vertices  $P\in X_1\cap Y$, $Q\in X_2\cap Y$ and $R\in X_1\cap X_2$. Since $Y$ is \CAT\ subcomplex and $P, Q\in Y$, we have ${\rm I}(P,Q)\subseteq Y$ and  therefore also $m(P,Q,R)\in Y$.
For the same reason $P, R\in X_1$ gives $m(P,Q,R)\in X_1$ and $Q, R\in X_2$ implies $m(P,Q,R)\in X_2$. We have shown that $m(P,Q,R)\in X_1\cap X_2\cap Y$. This finishes the proof.
\end{pf}

We need the following crucial definitions.
\begin{Definition}
\begin{enumerate}
\item[$(i)$] A simplicial action on a \CAT\ cube complex is called strongly simplicial if the stabilizer group of any cube fixes that cube pointwise.
\item[$(ii)$] A group $G$ is said to have property $\F^{s}\mathcal{C}$ if any strongly simplicial action of $G$ on any \CAT\ cube complex has a fixed point. 
\end{enumerate}
\end{Definition} 

Using Helly's Theorem for \CAT\ cube complexes  we obtain the following group theoretic condition for the fixed point property. 
\begin{proposition}
\label{HellyGroup1}
Let $G=\langle\left\{g_{1},\ldots, g_{k}\right\}\rangle$ be a finitely generated group and $X$ a  \CAT\ cube complex. If
$\Phi: G\rightarrow \isom$ is a strongly simplicial action such that the fixed point sets $\Fix(\langle g_{i}, g_{j}\rangle)$ are non-empty for all $i, j=1,\ldots, k$, then the fixed point set ${\rm Fix}(G)$ is non-empty as well.
\end{proposition}
\begin{pf}
Since the action is strongly simplicial, the fixed point sets of the generators are \CAT\ subcomplexes. By assumption we have that $\Fix(\langle g_{i}, g_{j}\rangle)=\Fix(g_i)\cap\Fix(g_j)$ is non-empty for all $i, j=1,\ldots, k$. It hence follows from Helly's Theorem for \CAT\ cube complexes  that $\Fix(g_1)\cap\ldots\cap\Fix(g_k)={\rm Fix}(G)$ is non-empty.
\end{pf}
\subsection{Fixed Point Criterion}
Combining the above proposition with Proposition \ref{boundedOrbit}, we immedia-\\tely find
\begin{CellularProposition} 
\label{CellularFPC}
Let $G$ be a group and $Y$ a finite generating set of $G$. If each pair of elements in $Y$ generates a finite subgroup, then $G$ has property $\F^{s}\mathcal{C}$.
\end{CellularProposition}
\begin{pf}
By Proposition \ref{boundedOrbit}, $\Fix(\langle g_{i}, g_{j}\rangle)$ is non-empty for all $i, j\in\left\{1,\ldots, k\right\}$. Now apply Proposition~\ref{HellyGroup1}.
\end{pf}

\begin{example}
\label{example}
Let $I$ be a finite set. A symmetric matrix $M=(m_{ij})_{ij}$ indexed by $I\times I$, with entries in $\mathbb{N}\cup\left\{\infty\right\}$ is called  Coxeter matrix if the following two properties are satisfied: $m_{ij}\geq 2$ for all $i\neq j$, $i, j\in I$ and $m_{ii}=1$ for all $i\in I$. The corresponding  Coxeter group $W$ has generating set $I$ and relators $(ij)^{m_{ij}}$. 
Every Coxeter group with a Coxeter matrix $M$ such that $m_{ij}<\infty$ for all $i, j\in I$ satisfies the Fixed Point Criterion and hence has property $\F^{s}\mathcal{C}$. But one can construct for every infinite Coxeter group a fixed point free simplicial action on a \CAT\ cube complex, see \cite{Coxeter}. Therefore an infinite  Coxeter group with finite entries in the Coxeter matrix has property $\F^{s}\mathcal{C}$ but doesn't have property $\F\mathcal{C}$.
\end{example}

\section{Generation \Aut\ and \SAut\ by finite subgroups}
We begin with the definition of the automorphism group of the free group of rank~$n$.
Let $F_{n}$ be the free group of rank $n$ with a fixed basis $X:=\left\{x_{1}, \ldots, x_{n}\right\}$. We denote by \Aut\ the automorphism group of $F_n$ and by \SAut\ the unique subgroup of index two in \Aut. More precisely, the abelianization map $F_{n}\twoheadrightarrow \Z^{n}$ gives a natural surjection $\pi:\Autn\twoheadrightarrow\GL_{n}(\Z).$
The group \SAut\ is equal to the preimage of $\SL_{n}(\Z)$ under this map. 

\subsection{A generating set of $\mathbf{Aut}\boldsymbol{(F_n)}$ and $\mathbf{SAut}\boldsymbol{(F_n)}$}

It was proven in \cite{VargheseFAd} that \Aut\ and \SAut\ have generating sets such that each pair of its elements generates a finite subgroup. 

\begin{proposition}(\cite{VargheseFAd})
\label{GenAut}
\begin{enumerate}
 \item[$(i)$] Let $n\geq 3$. The group \Aut\ has a generating set $Y_1$ such that any subgroup generated by $\left\{\alpha, \beta\right\}\subseteq Y_1$ is finite.
 \item[$(ii)$] Let $n\geq 4$. The group \SAut\ has a generating set $Y_2$ such that any subgroup generated by $\left\{\alpha, \beta\right\}\subseteq Y_2$ is finite.
\end{enumerate}
\end{proposition}

\section[Proof of Theorem A]{Proof of Theorem A}

Now we have all the ingredients to prove Theorem A. We show that the Fixed Point Criterion is satisfied by \Aut\ and \SAut, and therefore these groups have property $\F^{s}\mathcal{C}$.
\begin{NewTheoremA}
$ $
\label{AutFC}
\begin{enumerate}
\item[$(i)$] For $n\geq 3$ the group \Aut\ has property $\F^{s}\mathcal{C}$. 
\item[$(ii)$] For $n\geq 4$ the group \SAut\  has property $\F^{s}\mathcal{C}$. 
\end{enumerate}
\end{NewTheoremA}

\begin{pf}
The generating sets $Y_1$ and $Y_2$ in Proposition \ref{GenAut} show that \Aut\ for $n\geq 3$ and \SAut\ for $n\geq 4$ satisfy the Fixed Point Criterion. Therefore, these groups have property $\F^{s}\mathcal{C}$.
\end{pf}

\section{Triviality for actions of \SAut\ on \CAT\ cube complexes}
The aim of this section is to show that  \SAut\ cannot act non-trivially on an $m$-ary \CAT\ cube complex for $m<n$. By definition, an $m$-ary cube complex is a cube complex where every vertex has at most $m$ neighbours. 

For the proof we need the following variant of a result by \textsc{Bridson} and \textsc{Vogtmann} \cite[3.1]{VogtmannSpheres}. For a detailed proof the reader is referred to \cite[1.13]{Diplomarbeit}.
\begin{proposition}
\label{faktorisiertSL}
Let $n\geq 3$, $G$ be a group and $\phi : \SAutn\rightarrow G$ a group homomorphism. 
 If there exists $\alpha\in{\rm Alt}(n)-\left\{{\rm id}_{F_n}\right\}$ with $\phi(\alpha)=1$, then $\phi$ is trivial.
\end{proposition}
We finish by proving
\begin{NewTheoremB}
Let $n\geq 4$ and $X$ be an $m$-ary \CAT\ cube complex. Let $\Phi:\SAutn\rightarrow{\rm Isom}(X)$
be a strongly simplicial action. If $m<n$, then $\Phi$ is trivial.
\end{NewTheoremB}
\begin{pf}
The group \SAut\ has property $\F^{s}\mathcal{C}$ by Theorem A, therefore $\Phi$ has a global fixed vertex $v\in X$. The group \SAut\ acts on the link of~$v$, i. e. the set of all neighbours of~$v$, via ${\rm Sym}(m)$. As \SAut\ is perfect, we even have
\[
\Phi_{{\rm Alt}(m)}:\SAutn\rightarrow{\rm Alt}(m).
\]
If $m<n$, the restriction of this map to ${\rm Alt}(n)$ cannot be injective, therefore  $\Phi_{{\rm Alt}(m)}$ is trivial by Proposition \ref{faktorisiertSL}. This shows that any neighbour of $v$ is in the fixed point set of the action, hence all vertices of $X$ are in the fixed point set. Thus \SAut\ acts trivially on $X$.
\end{pf}

\subsection*{Acknowledgements} The author would like to thank Yves Cornulier and Genevois Anthony for their comments concerning completeness in Proposition \ref{boundedOrbit} and the referee for many helpful comments.

\end{document}